\documentclass[11pt]{article}
\usepackage{caption2}
\usepackage{pdfsync}
\usepackage{amssymb}
\usepackage{color}
\usepackage{amsmath,amsthm}
\usepackage{graphicx}
\usepackage{empheq}
\usepackage{indentfirst}
\def\f{\frac}

\def\i1n{i=1,\cdots,n}
\def\j1n{j=1,\cdots,n}
\def\ij1n{i,j=1,\cdots,n}

\def\R{\mathbb R}

\newcommand{\beq}{\begin{equation*}}
\newcommand{\eeq}{\end{equation*}}
\newcommand{\be}{\begin{equation}}
\newcommand{\ee}{\end{equation}}
\newcommand{\tr}{^\text{tr}}

\newcommand{\cqfd}
{%
\mbox{}%
\nolinebreak%
\hfill%
\rule{2mm}{2mm}%
\medbreak%
\par%
}
 \newtheorem{thm}{Theorem}[section]
 
 \newtheorem{lem}{Lemma}[section]

 \newtheorem{rem}{Remark}[section]
\usepackage{color}
\usepackage{float}
\usepackage{graphicx}
\usepackage{indentfirst}
\usepackage{layout}

\def\R{\mathbb R}

\title{Optimal control of cell mass and maturity in a model of follicular ovulation}

\author{Fr\'{e}d\'{e}rique Cl\'{e}ment\thanks{INRIA Paris-Rocquencourt Centre.
Rocquencourt BP 105-78153 Le Chesnay Cedex, France.
E-mail: {\tt Frederique.Clement@inria.fr}.}
        \and Jean-Michel Coron\thanks{Institut universitaire de France and Universit\'{e} Pierre et
Marie Curie-Paris 6, UMR 7598 Laboratoire Jacques-Louis
Lions, 75005 Paris, France. E-mail: {\tt coron@ann.jussieu.fr}. JMC was partially
supported by the ERC advanced grant 266907 (CPDENL) of the 7th Research Framework Programme (FP7).}
\and Peipei Shang\thanks{INRIA Paris-Rocquencourt Centre. Universit\'{e}
Pierre et Marie Curie-Paris 6, UMR 7598 Laboratoire Jacques-Louis
Lions, 75005 Paris, France. E-mail: {\tt Peipei.Shang@inria.fr}. PS
was supported by the INRIA large scale initiative action \href{https://www.rocq.inria.fr/sisyphe/reglo/regate.html}{REGATE} (REgulation of
the GonAdoTropE axis).}}
\date{}
\begin{document}

\maketitle

\begin{abstract}
In this paper, we study optimal control problems associated with a scalar hyperbolic conservation law modeling the development of ovarian follicles. 
Changes in the age and maturity of follicular cells are described by a 2D conservation law, where the control terms act on the velocities. 
The control problem consists in optimizing the follicular cell resources so that the follicular maturity reaches a maximal value in fixed time. 
Formulating the optimal control problem within a hybrid framework, we prove necessary optimality conditions in the form of Hybrid Maximum Principle. 
Then we derive the optimal strategy and show that there exists at least one optimal bang-bang control with one single switching time. 
\end{abstract}


{\bf 2000 MR Subject Classification:}\quad
35L65, 49J20, 92B05.

\section{Introduction}

This work is motivated by natural control problems arising in reproductive physiology. The development of ovarian follicles is a crucial process for reproduction in mammals, as its biological meaning is to free fertilizable oocyte(s) at the time of ovulation. 
During each ovarian cycle, numerous follicles are in competition for their survival. Few follicles reach an ovulatory size, since most of them undergo a degeneration process, known as atresia (see for instance \cite{MH}). The follicular cell population consists of proliferating, differentiated and apoptotic cells, and the fate of a follicle is determined by the changes occurring in its cell population in response to an hormonal control originating from the pituitary gland.

A mathematical model, using both multiscale modeling and control theory concepts, has been designed to describe the follicle selection process on a cellular basis (see \cite{FC05}).
The cell population dynamics is ruled by a conservation law, which describes the changes in the distribution of cell age and maturity.

Cells are characterized by their position within or outside the cell cycle and
by their sensitivity to the follicle stimulating hormone (FSH). This leads one to distinguish 3 cellular phases. Phase 1 and 2 correspond to the proliferation phases and Phase 3 corresponds to the differentiation phase, after the cells have exited the cell cycle.

The cell population in a follicle $f$ is represented by cell density functions $\rho^f_{j,k}(t,x,y)$
defined on each cellular phase $Q^f_{j,k}$, where $j=1,2,3$ denotes
Phase 1, Phase 2 and Phase 3, $k=1,2,\cdots$ denotes the number of the successive cell cycles (see figure~\ref{chart}).
The cell density functions satisfy the following conservation laws:
\be
\label{ocl}
\frac{\partial\rho^f_{j,k}}{\partial t}+\frac{\partial(g_f(u_f)\rho^f_{j,k})}{\partial x}+\frac{\partial(h_f(y, u_f)\rho^f_{j,k})}{\partial y}=-\lambda(y, U)\rho^f_{j,k}\,
\text{ in }\, {Q^f_{j,k}},
\ee
where $Q^f_{j,k}=\Omega^f_{j,k}\times[0,T]$, with
\begin{align*}
&\Omega^f_{1,k}=[(k-1)a_2,(k-1)a_2+a_1]\times[0,y_s], \\
&\Omega^f_{2,k}=[(k-1)a_2+a_1,ka_2]\times[0,y_s],\\
&\Omega^f_{3,k}=[(k-1)a_2,ka_2]\times[y_s,y_m].
\end{align*}
\vspace{10mm}
\setlength{\unitlength}{0.085in}
\begin{figure}[htbp]
\begin{picture}(-10,10)
\put(10,1){\vector(1,0){50}} \put(10,1){\vector(0,1){15}}

\put(10,13){\line(1,0){45}} \put(10,7){\line(1,0){45}}
\put(15,1){\line(0,1){5.9}} \put(20,1){\line(0,1){12}}
\put(25,1){\line(0,1){5.9}} \put(30,1){\line(0,1){12}}
\put(38,1){\line(0,1){12}}  \put(55,1){\line(0,1){12}}

\put(33,3){\makebox(2,1)[l]{$\cdots$}}
\put(33,8){\makebox(2,1)[l]{$\cdots$}}

 \put(43,1){\line(0,1){5.9}}
\put(48,1){\line(0,1){12}}

\put(7.5,6){\makebox(2,1)[l]{$y_s$}}
\put(7,12){\makebox(2,1)[l]{$y_m$}}
\put(9.5,17){\makebox(2,1)[l]{$y$}}

\put(61,0.4){\makebox(2,1)[l]{$x$}}
\put(9,-1){\makebox(2,1)[l]{$0$}}
\put(14,-1){\makebox(2,1)[l]{\tiny$a_1$}}
\put(19,-1){\makebox(2,1)[l]{\tiny$a_2$}}
\put(23,-1){\makebox(2,1)[l]{\tiny$a_2+a_1$}}
\put(29,-1){\makebox(2,1)[l]{\tiny$2a_2$}}
\put(34.5,-1){\makebox(2,1)[l]{\tiny $(k-1)a_2$}}
\put(40,-1){\makebox(2,1)[l]{\tiny $(k-1)a_2+a_1$}}
\put(48,-1){\makebox(2,1)[l]{\tiny $ka_2$}}

\put(11,3){\makebox(2,1)[l]{ $\rho^f_{1,1}$}}
\put(16.5,3){\makebox(2,1)[l]{$\rho^f_{2,1}$}}
\put(14,9.5){\makebox(2,1)[l]{$\rho^f_{3,1}$}}

\put(21,3){\makebox(2,1)[l]{ $\rho^f_{1,2}$}}
\put(26.5,3){\makebox(2,1)[l]{$\rho^f_{2,2}$}}
\put(25,9.5){\makebox(2,1)[l]{$\rho^f_{3,2}$}}

\put(39,3){\makebox(2,1)[l]{ $\rho^f_{1,k}$}}
\put(44.5,3){\makebox(2,1)[l]{$\rho^f_{2,k}$}}
\put(42,9.5){\makebox(2,1)[l]{$\rho^f_{3,k}$}}

\put(50.5,3){\makebox(2,1)[l]{$\cdots$}}
\put(50.5,8){\makebox(2,1)[l]{$\cdots$}}

\put(19,7.5){\vector(1,1){1.7}} \put(19,12.5){\vector(1,-1){1.7}}
\put(29,7.5){\vector(1,1){1.7}} \put(29,12.5){\vector(1,-1){1.7}}
\put(37,7.5){\vector(1,1){1.7}} \put(37,12.5){\vector(1,-1){1.7}}
\put(47,7.5){\vector(1,1){1.7}} \put(47,12.5){\vector(1,-1){1.7}}

\put(12,6){\vector(1,1){1.7}}
\put(14,4){\vector(1,1){1.7}}
\put(19,4){\vector(1,0){1.5}}
\put(22,6){\vector(1,1){1.7}}
\put(24,4){\vector(1,1){1.7}}
\put(29,4){\vector(1,0){1.5}}
\put(40,6){\vector(1,1){1.7}}
\put(37,4){\vector(1,0){1.5}}
\put(42,4){\vector(1,1){1.7}}
\put(47,4){\vector(1,0){1.5}}
\end{picture}
\vspace{10mm}
\caption{Cellular phases on the age-maturity plane
for each follicle $f$. The domain consists of the sequence of
$k=1,2,\cdots$ cell cycles. Variable $x$ denotes the age of the cell and variable $y$ denotes its maturity. Number $y_s$ is the threshold value at which cell cycle exit occurs and  $y_m$ is the maximal maturity.
The top of the domain corresponds to the differentiation phase and the bottom to the proliferation phase.}
\label{chart}
\end{figure}
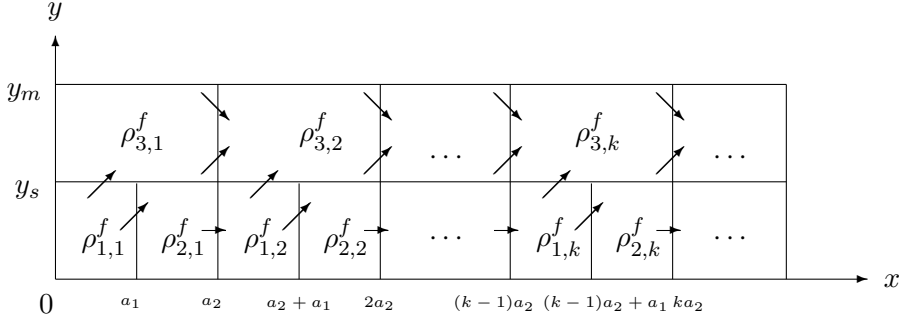

Let us define
\be\label{moment}
M_f(t):=\sum_{j=1}^3\sum_{k=1}^N\int_0^{+\infty}\int_0^{+\infty}y\, \rho^f_{j,k}(t,x,y)\,dx\, dy
\ee
as the maturity on the follicle scale, and
\be\label{moment2}
M(t):=\sum_f M_f(t)
\ee
as the maturity on the ovarian scale.

The velocities of aging $g_f$ and maturation $h_f$ as well as the loss term $\lambda$ depends on the mean maturity of the follicle $f$
through a local control $u_f(t,M_f,M)$ which represents intrafollicular bioavailable FSH level and the mean maturity of all the follicles through a global
control $U(t,M)$ which can be interpreted as the FSH plasma level.
One can refer to \cite{FC07, FC05, Shang2} for more details on the model.

The aging velocity controls the duration of the cell division cycle. Once the cell age has reached a critical age, the mitosis event is triggered and a mother cell gives birth to two daughter cells.
The two
daughter cells enter a new cell cycle, which results in
a local doubling of the flux. Hence, there are local singularities in the subpart of the domain where $y\leqslant y_s$, that correspond to the flux doubling due to the successive mitosis events. 
 The maturation velocity controls the time needed to reach a threshold maturity $y_s$, when the cell exits the division cycle definitively.  After  the exit time, the cell is no more able to contribute to the increase in the follicular cell mass.

Ovulation is triggered when the ovarian maturity reaches a threshold value $M_s$. The stopping time $T_s$ is defined as
\be
\label{ot} T_s:=\min\left\{T\ | M(T)=M_s\right\},
\ee
and corresponds on the biological ground to the triggering of a massive secretion of the hypothalamic gonadotropin releasing hormone (GnRH).

As a whole, system \eqref{ocl}-\eqref{moment2} combined with stopping condition \eqref{ot} defines a multiscale reachability problem. It can be associated to an optimal control problem that consists in minimizing $T_s$ for a given target maturity $M_s$.

Some related control problems have already been investigated on a mathematical ground. In \cite{FC07}, the authors studied the characteristics associated with a follicle  as an open-loop control problem. They described the sets of microscopic initial conditions compatible with either ovulation or atresia in the framework of backwards reachable set theory. Since these sets were largely overlapping, their results illustrate the prominent impact of cell dynamics control in the model. In \cite{Michel}, the author focused on the issue of the selection process in a game theory approach, where one follicle plays against all the other ones. Whether the follicle becomes atretic (doomed) or ovulatory (saved) depends on the follicular cell mass reached at the time when all cells stop proliferating.

The aim of this paper is to investigate whether there exists an optimal way
for a follicle to reach ovulation.
On the one hand, the follicle can benefit from a
strong and quick enlargement of its cell population. On the other hand, this enlargement 
occurs at the expense of the maturation of individual cells. This compromise
was instanced here as a problem of composition of velocities.
A concept central to the understanding of these entangled processes is that of the management of follicular cell resources.
There is indeed a finely tuned balance between the production of new cells through proliferation, that increases the whole cell mass, and the maturation of cells, that increases their contribution to hormone secretion.

The controllability of nonlinear
hyperbolic equations (or systems) have been widely studied for a long time; for the 1D case,
see, for instance \cite{Coron,
CGWang, CW2012, 2003-Gugat-SICON, 2009-Gugat-Leugering-AIHP, LiBook1, LiRao, LRW2010, Wang} for smooth solutions and
\cite{Ancona98, 2002-Bressan-Coclite-SICON, Glass, Horsin} for bounded variation entropic solutions. In particular,
\cite{CoronBook} provides a comprehensive survey of controllability
of partial differential equations including nonlinear hyperbolic systems. As far as optimal control problems for hyperbolic systems are concerned, one can refer to \cite{Gugat2008, GHKL, GHS, SW}.
However, most of these monographs study the case where the controls are either applied inside the
domain or on the boundary.
Our control problem is quite different from the problems already studied in the literature, since the control terms appear in the flux.
To solve the problem, we make use both of analytical methods based on Hybrid Maximum Principle (HMP) and numerical computations.

The paper is organized as follows. In section \ref{pre}, we set the optimal control problem, together with our assumptions, and we enunciate the main result. In section \ref{Dirac}, we give necessary optimality conditions from HMP in the case where Dirac masses are used as a rough approximation of the density. 
An alternative sketch of the proof based on an approximation method is given in appendix.
Using the optimality conditions, 
we show that for finite Dirac masses, every measurable optimal control is a bang-bang control with one single switching time. 
In addition to the theoretical results, we give some numerical illustrations. In section \ref{Pde}, we go back to the original PDE formulation of the model, and we show that there exists at least one optimal bang-bang control with one single switching time.

\section{Problem statement and introductory results}\label{pre}
\subsection{Simplifications with respect to the original model}\label{pre2}
To make the initial problem tractable, we have made several simplifications on the model dynamics.
\begin{align*}
S_1. &\text{ We consider only one developing follicle},\ \text{i.e.}\ f=1;\\
S_2. &\text{ There is no loss term anymore},\ \text{i.e.}\ \lambda=0;\\
S_3. &\text{ The age velocity is uncontrolled},\ \text{i.e.}\  g_f \equiv 1;\\
S_4. &\text{ The cell division is represented by a new gain term},\ \text{i.e.}\  c(y) \text{ defined by }\eqref{scp};\\
S_5. &\text{ The target maturity $M_s$ can always be reached in finite time.}
\end{align*}

($S_1$) means that, in this problem, we are specially interested in the control of the follicular cell resources for each follicle,
in the sense that we 
ignore the influence of the
other growing follicles.
The goal is to find the optimal balance between the production of new cells and the maturation of cells.

In ($S_2$), we neglect the cell death, which is quite natural when considering only ovulatory trajectories, while, in ($S_3$), we consider that the cell age evolves as time. Moreover, the cell division process is distributed over ages with ($S_4$), so that there is a new gain term in the model instead of the former mitosis transfer condition. 

Even if it is simplified, the problem studied here still captures the essential
 question of the compromise between proliferation and differentiation that
 characterizes terminal follicular development.
 A relatively high aging velocity tends to favor cell mass production, 
while a relatively high maturation velocity tends to favor an increase in the average cell maturity.

As shown in section \ref{sectionduality}, assumptions
($S_2$) and ($S_5$) allow us to replace a minimal time criterion by a criterion that consists in maximizing the final maturity.
Hence, from the initial, minimal time criterion, we have shifted, for sake of technical simplicity, to an equivalent problem where the final time is fixed and the optimality criterion is the follicular maturity at final time.
On the biological ground, this means that for any chosen final time $t_1$,
the resulting maturity at final time $M_f (t_1)$ can be chosen in turn as a maturity target which would be reached in minimal time at time $t_1$.
It can be noticed that in the initial problem \eqref{ot}, there might be no optimal solution without assumption ($S_5$), if the target maturity is higher than the maximal asymptotic maturity.

\subsection{Optimal control problem}\label{pre1}
Under these assumptions,
we arrived to consider the following conservation law on a fixed time horizon:
\begin{gather}
\label{seqCauchy}
\left\{
\begin{array}{l}
\rho_t +\rho_x +((a(y)+b(y)u)\rho)_y=c(y)\rho,\quad t\in(t_0,t_1),\,
x> 0,\, y> 0,
\\
\rho(t,0,y)=\rho(t,x,0)= 0,\quad  t\in (t_0,t_1),\, x> 0,\, y> 0,
\\
\rho(0,x,y)=\rho_0(x,y),\quad  x> 0,\, y> 0,
\end{array}
\right.
\end{gather}
where
\begin{align}\label{ab}
a(y):=-y^2,\quad b(y):=c_1 y+c_2,
\end{align}
and
\be
\label{scp}
c(y):=\begin{cases}\ c_s,\quad\text{if}\ y\in[0,y_s),\\
\ 0,\quad\text{if}\ y\in[y_s,\infty),
\end{cases}
\ee
with $y_s$, $c_s$, $c_1$ and $c_2$ being given strictly positive constants. We assume that
\begin{gather}
\label{assumptionys}
\f{y^2_s}{c_1y_s+c_2}<1.
\end{gather}
Let us denote by $w$ a positive constant such that
\be
\label{w}
w\in (\f{y^2_s}{c_1y_s+c_2},1).
\ee
 From \eqref{ab} and \eqref{w}, we have
\be
\label{h(ys)>0}
a(y)+b(y)u>0, \, \forall y \in [0,y_s], \, \forall u\in [w,1].
\ee
Throughout this paper the control $u$ is assumed to satisfy the constraint
\begin{gather}
\label{uin[w,1]}
u\in [w,1].
\end{gather}
The left constraint $w$
in \eqref{uin[w,1]} ensures that the maturation velocity is always positive in the proliferation
phase.
The right constraint in \eqref{uin[w,1]} is natural since FSH plasma levels are bounded.
The maximal bound can be scaled to $1$ for sake of restricting the number of parameters in the model.

By \eqref{uin[w,1]}, 
there is a maximal asymptotic maturity $\bar y$ on the cell scale,
 i.e. the positive root $y$ of $a(y)+ b(y)u=0$ with control $u=1$. From  \eqref{ab},
we have
\be
\label{bary}
\bar
y=\displaystyle\f{c_1+\sqrt{c^2_1+4c_2}}{2}.
\ee

Let $u\in L^\infty((t_0,t_1);
[w,1])$. Let us define the map
\begin{equation*}
\begin{array}{ccc}
\Psi: [t_0,t_1]\times [0,y_s]\times L^\infty ((t_0,t_1);[w,1])
&
\rightarrow
&
[0,\bar y]
\\
(t,y_0,u)
&
\mapsto
&
\Psi (t,y_0,u)
\end{array}
\end{equation*}
by requiring
\begin{equation}
\label{defPhi}
\left\{
\begin{array}{l}
\displaystyle
\frac{\partial \Psi}{\partial t} (t,y_0,u)= a(\Psi (t,y_0,u))+ b(\Psi (t,y_0,u))u(t),
\\
\Psi (t_0,y_0,u)=y_0.
\end{array}
\right.
\end{equation}
Let us now define the exit time $\hat
t_0$ as
\be\label{hat0}
\Psi(\hat t_0,0,w)=y_s.
 \ee
Let us point out that, by \eqref{h(ys)>0}, there exists one and only one  $\hat t_0$
satisfying \eqref{hat0}. Note that it is not guaranteed that the exit time $\hat t_0$ occurs
before the final time $t_1$, so that we may have $\hat t_0>t_1$. 
When $t>\hat t_0$, all the cells are in Phase 3, i.e. their maturity
is larger than the threshold $y_s$. After time $\hat t_0$
the mass will not increase any more due to  \eqref{scp}.
The maximal cell mass that can be reached at $\hat t_0$ is obtained 
when applying $u=w$ from the initial time.

For any admissible control $u\in L^{\infty}((t_0,t_1);[w,1])$, we define the
cost function
\be
\label{sc1defJ}
J(u):=-\int_0^{+\infty}\int_0^{+\infty}y\, d\rho(t_1,x,y),
\ee
and we want to study the following optimal control problem:
\be\label{sc1minJ}
\text{minimize } J(u) \text{ for } u\in L^{\infty}((t_0,t_1);[w,1]).
\ee

A similar minimal time problem was investigated in another ODE framework \cite{FC08}, where the proliferating and differentiated cells were respectively pooled in a proliferating and a differentiated compartment. 
The author proved by Pontryagin Maximum Principle (PMP) that the optimal strategy is a bang-bang control, which consists in applying permanently the minimal apoptosis rate and in switching once the cell cycle exit rate from its minimal bound to its maximal one.
In contrast, due to the fact that $c$ is discontinuous, we cannot apply
PMP directly here. The idea is to
first consider optimal control problems for Dirac masses (see section
\ref{Dirac}), and then to pass to the limit to get optimal control results for
the PDE case (see section
\ref{Pde}). 

For ``discontinuous'' optimal control problems of finite dimension,
one cannot derive necessary optimality conditions by applying directly the
standard apparatus of the theory of extremal problems \cite{BP, AV, LR}.
The first problem where the cost function was an integral functional with discontinuous
integrand was dealt in \cite{AVA}.
Later, in \cite {S}, the author studied the case of a more general functional that includes
both the discontinuous characteristic function and continuous terms.
There, the author used approximation methods to prove necessary optimality
conditions in the form of PMP. One of the difficulties of our problem is that both the integrand of the cost function and the dynamics are discontinuous.

However, 
our problem can be classified as a hybrid optimal control problem,
since the problem has a discontinuous dynamics ruled by a partition of the state space.
One of the most important results in the study of such problems is the HMP
proved in \cite{Piccoli2, Piccoli1, Sussmann2, Sussmann1}.
There, the authors followed the standard line of the full procedure for the direct proof of PMP,
based on the introduction of a special class of control variations, and the computation of the increments of the cost and all constraints.
In \cite{Dmitruk}, the authors formulated the hybrid problem as a classical optimal control problem.
They then proved the HMP using the classical PMP.
Later, in \cite{TH}, the authors regularized the hybrid problems to standard smooth optimal control problems,
to which they can apply the usual PMP. They also derived jump conditions appropriate to our problem.

The main result of this paper is the following theorem. 
\begin{thm}\label{Theorempde}
Let us assume that
\begin{gather}
\label{t1>hatt0}
t_1>\hat t_0,
\\
2y_s-c_1>0\quad\text{and}\quad c_s> \f{a(y_s)+b(y_s)}{y_s}.
\label{Dsuse}
\end{gather}
Then, among all admissible controls $u\in L^\infty((t_0,t_1);
[w,1])$, there exists an optimal control $u_*$ for the minimization problem \eqref{sc1minJ} such that
\begin{gather}
\label{bang-bang-effectif-1}
\exists\, t_*\in [t_0,t_1] \text{  such that } u_*= w \text{ in } (t_0,t_*)\ \text {
and }\ u_*=1 \text{ in } (t_*,t_1).
\end{gather}
\end{thm}

\begin{rem}
From the mathematical viewpoint, assumptions \eqref{t1>hatt0} and \eqref{Dsuse}
arise naturally from the computations (see section \ref{one}).
Condition \eqref{t1>hatt0} means that we consider a target time large enough so that all the cells have gone to the differentiation phase.
Condition \eqref{Dsuse} gives specific relations between the proliferation rate and the parameters of the maturation velocity.
Together, these relations are related to the transit time within the proliferation phase.
\end{rem}
\begin{rem}
In our case, the dynamics of $\rho$ is essentially one-dimensional,
since there is a transport
with constant velocity along variable $x$
and we have just to deal with variable $y$.
Hence our results can be generalized to $n$-spatial dimensional problem like
\be
\rho_t+m\cdot \nabla_x \rho +(h(y,u)\rho)_y=c(y)\rho,
\ee
where $m$ is a constant vector.
Generalization to $n$-spatial dimensional dynamics with both velocities controlled should also be feasible.
\end{rem}

\subsection{Solution to Cauchy problem \eqref{seqCauchy}}
\label{introre}

In this section, we give the definition of a (weak) solution to Cauchy problem \eqref{seqCauchy}.

\begin{figure}[htbp!]
\includegraphics[width=1.7in]{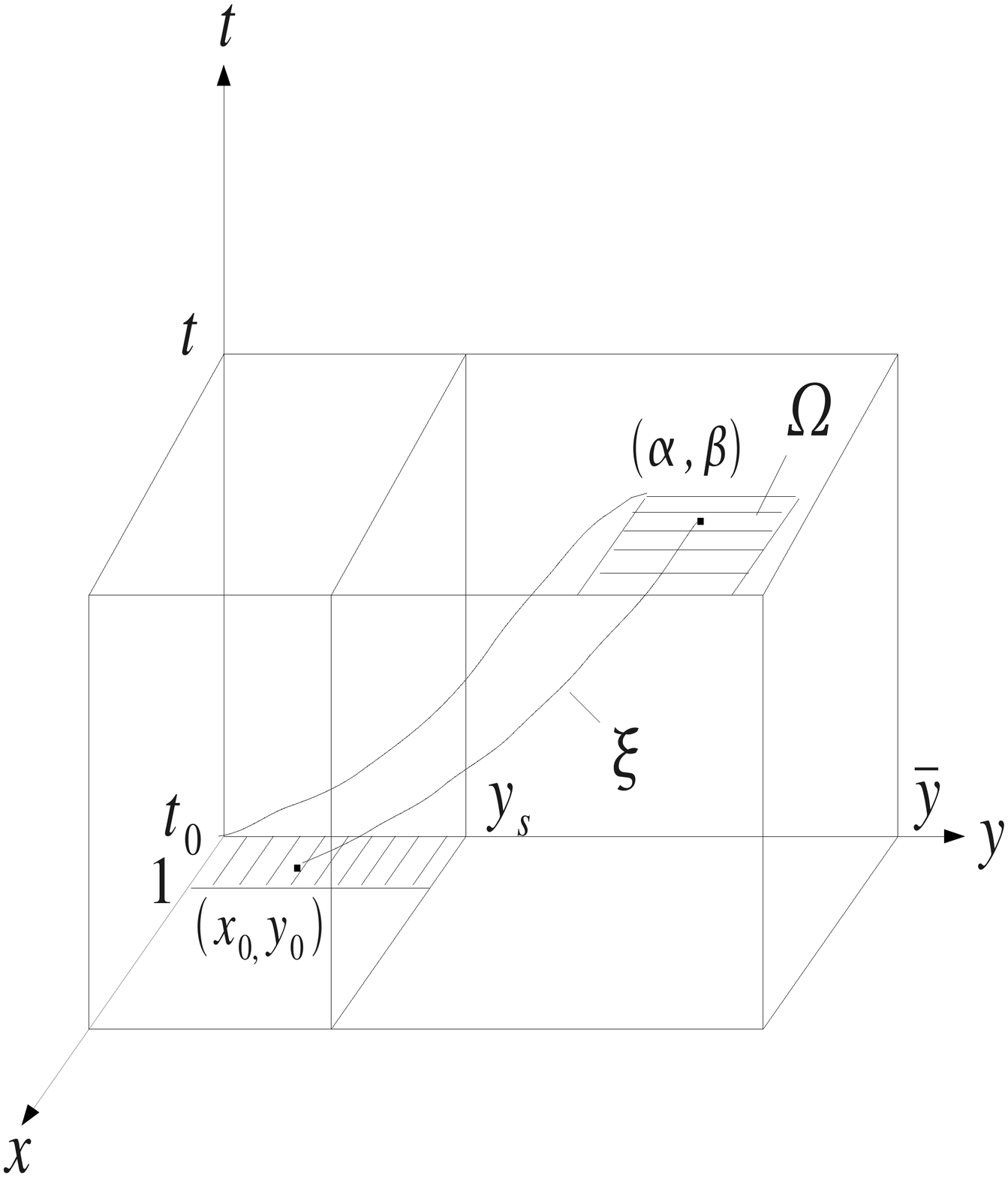}
\centering
\caption{Construction of a weak solution by backward tracking of the characteristics. Variables $x$ and $y$ respectively denote the age
and maturity, $y_s$ is the threshold maturity and $\bar y$ is the maximal asymptotic maturity.
The initial mass concentrates in the shaded area $[0,1]\times[0,y_s]$. The curve
$\xi=(x,y)$ is the characteristic curve passing through $(t,\alpha,\beta)$ that intersects the initial plane $t=t_0$ at $(t_0,x_0,y_0)$.}
\label{charac}
\end{figure}

Let $\rho_0$ be a Borel measure on $\R\times\R$  such that
\begin{gather}
\rho_0\geqslant 0,
\label{rho0>0}
\\
\text{and the support of $\rho_0$ is included in $[0,1]\times [0,y_s]$}.
\label{supportrho0}
\end{gather}
Let $K:=[0,t_1-t_0+1]\times[0,\bar y]$. Let $M(K)$ be the set  of Borel measures on $K$, i.e.  the set of continuous linear maps from $C^0(K)$ into $\R$. The solution to Cauchy problem \eqref{seqCauchy} is the function $\rho:\, [t_0,t_1]\rightarrow M(K)$ such that, for every $\varphi\in C^0(K)$,
 \be \label{solution}
\iint_{K} \varphi(\alpha,\beta)d\rho(t,\alpha,\beta)
=\iint_K \varphi(x_0+t-t_0,\Psi(t,y_0,u))\, e^{\int_{t_0}^{t}c(\Psi(s,y_0,u))ds
}\, d\rho_0(x_0,y_0).
\ee
We take expression \eqref{solution} as a definition. This expression is also justified by the fact that if $\rho_0$ is a $L^\infty$ function, one recovers the usual notion of weak solutions to Cauchy problem \eqref{seqCauchy} studied in \cite{CoronBook, CKW, Shang2, SW},
as well as by the characteristics method used to solve hyperbolic equations (see figure~\ref{charac}).

\begin{rem}
From \eqref{solution},
if $\rho_0$ is a positive Borel measure, then solution $\rho$ 
is also a positive Borel measure.
If $\rho_0\in C^0$ or if $\rho_0$ is Lipschitz continuous, then $\rho\in C^0$ or $\rho$ is Lipschitz continuous, but 
if $\rho_0\in C^1$, it may happen that $\rho$ is not in $C^1$ due to the fact that $c(y)$ is discontinuous.
\end{rem}

\subsection{Minimal time versus maximal maturity}
\label{sectionduality}
In this section, we show that the two optimal control problems enunciate either as: ``minimize the time to achieve a given maturity'' or
``achieve a maximal maturity at a given time'' are equivalent when $S_2$ and $S_5$ hold.
The threshold target maturity $\bar M_s$ in $S_5$ can be computed from the maximal cell mass combined with the maximal asymptotic
maturity $\bar y$ when applying $u=w$ from the initial time until $\hat t_0$ and $u=1$ thereafter, 
so that $S_5$ can be formulated as:
\be\label{asymaturityfollicle}
M_s< \bar M_s:=\bar y\int_0^{\infty}\int_0^{\infty}d\rho(\hat t_0,x,y).
\ee

Let $\rho_0$ be a nonzero Borel measure on $\R\times\R$
satisfying \eqref{rho0>0} and \eqref{supportrho0}. Let us denote by $M^u(t)$ the maturity at time $t$ for the control
$u\in L^{\infty}((t_0,t_1);[w,1])$ (and the initial data $\rho_0$).

\text{A.} For fixed target time $t_1$, suppose that the maximum of the
maturity $M$
is achieved with an optimal control
$u\in L^{\infty}((t_0,t_1);[w,1])$
\be
M^{u}(t_1)= M.
\ee
Then we conclude that for this fixed $M$, the minimal time needed to
reach $M$ is $t_1$ with the same control $u$. We prove it by
contradiction. We assume that there exists another control $\tilde u \in L^{\infty}((t_0,\tilde t_1);[w,1])$ such that
\be
M^{\tilde u}(\tilde t_1)=M,\quad \tilde t_1<t_1.
\ee
We extend $\tilde u $ to $[t_0,t_1]$ by requiring $\tilde u = 1$ in $(\tilde t_1,t_1]$.
Let us prove that
\begin{gather}
\label{Mincreasing}
t\in [\tilde t_1,t_1] \rightarrow M^{\tilde u}(t)
\text{ is strictly increasing.}
\end{gather}
Let $\tilde \rho :[t_0,t_1]\rightarrow M(K)$ be the solution to Cauchy problem \eqref{seqCauchy} (see section~\ref{introre}).
Note that $a(y)+b(y)>0$ for every $y \in [0, \bar y)$ and that, for every $t\in [t_0,t_1]$,
  the support of $\tilde \rho(t)$ is included in $[0,t_1-t_0+1]\times [0,y_s)$. Together with \eqref{solution} for $\rho=\tilde \rho$
  and $\varphi(\alpha,\beta)=\beta$, this proves
\eqref{Mincreasing}. From \eqref{Mincreasing} it follows that
\be
M^{\tilde u}(t_1)>M^{\tilde
u}(\tilde t_1)=M,
\ee
which is a contradiction with the optimality of $u$.

\text{B.} For any fixed target maturity $M$, suppose that the minimal time needed to reach
$M$ is $t_1$ with control $u\in L^{\infty}((t_0, t_1);[w,1])$.
Then we conclude that for this fixed target time
$t_1$, the maximal maturity at time $t_1$ is $M$ with the same control $u$.
We prove it again by contradiction. We assume that there exists another control $\tilde
u\in L^{\infty}((t_0, t_1);[w,1])$ such that
\be
M^{\tilde u}(t_1)>M.
\ee
Then by the continuity of
$M^{\tilde u}(t)$ with respect to time $t$, there
exists a time $\tilde t_1<t_1$ such that
\be
M^{\tilde u}(\tilde
t_1)=M,
\ee
which is a contradiction with the minimal property of $t_1$.
This concludes the proof of the equivalence between the two optimal control problems.
\cqfd
\section{Results on optimal control for finite Dirac masses}\label{Dirac}
In this section, we give results on the optimal control problem \eqref{sc1minJ} when the initial data
$\rho_0\geqslant 0$ is a linear combination of a finite number of Dirac
masses. For $(\alpha,\beta)\tr \in\R^2 $, we denote by $\delta_{\alpha,\beta}$ the Dirac mass at $(\alpha,\beta)\tr$.
We assume that, for some positive integer $N$, there exist $N$ elements $((x^{k0}_1,x^{k0}_2))_{k\in\left\{1,\ldots,N\right\}}$
of $[0,1]\times[0,y_s]$  and $N$ strictly positive real numbers $(x^{k0}_3)_{k\in\left\{1,\ldots,N\right\}}$ such that
\be\label{initialD}
\rho_0:=\sum_{k=1}^N x^{k0}_3\delta_{x^{k0}_1,x^{k0}_2}.
\ee

First, we formulate our problem within a hybrid framework.
Let us denote by $X_{\alpha}$ and $X_{\beta}$ two disjoint and open subsets of $\R^{3}$, where
\begin{align*}
X_{\alpha}:=\Big\{(x_1,x_2,x_3)\in \R^{3}\, |x_2< y_s\Big\},\\
X_{\beta}:=\Big\{(x_1,x_2,x_3)\in \R^{3}\, |x_2> y_s\Big\}.
\end{align*}
The boundary between the two domains $X_{\alpha}(t)$ and $X_{\beta}(t)$ can be written as
\beq
\Big\{(x_1,x_2,x_3)\in \R^{3}\, |F(x)=0\Big\},
\eeq
where 
\be
F(x):=x_2-y_s.
\ee

We consider the following Cauchy problem:
\be\label{PN}
\begin{cases}
\dot x^k=f(x^k,u),\quad u\in L^{\infty}((t_0,t_1);[w,1]),\quad t\in[t_0,t_1],\\
x^k(t_0)=x^{k0},
\end{cases}
\ee
where
\be\label{vf}
f(x^k,u):=\begin{cases}
f_{\alpha}(x^k,u),\quad x^k\in X_{\alpha},\\
f_{\beta}(x^k,u),\quad x^k\in X_{\beta},
\end{cases}\quad
x^k=\begin{pmatrix}
x^k_1
\\
x^k_2
\\
x^k_3
\end{pmatrix}
,\,
x^{k0}=\begin{pmatrix}
x^{k0}_1
\\
x^{k0}_2
\\
x^{k0}_3
\end{pmatrix},
\ee
with
\begin{align*}
f_{\alpha}(x^k,u):=\left(
\begin{array}{c}
 1  \\
 a(x^k_2)+b(x^k_2)\, u \\
  c_s x^k_3
\end{array}
\right),
\;
f_{\beta}(x^k,u):=\left(
\begin{array}{c}
 1  \\
 a(x^k_2)+b(x^k_2)\, u \\
0
\end{array}
\right).
\end{align*}
It is easy to check that the maximal solution to Cauchy problem
 \eqref{PN}  is defined on
$[t_0,t_1]$.

One can also easily check that the solution to Cauchy problem \eqref{seqCauchy},
as defined in section~\ref{introre},  is
\be
\rho(t)=\sum_{k=1}^N x^{k}_3(t)\delta_{x^{k}_1(t),x^{k}_2(t)}.
\ee
The cost function $J$ defined in \eqref{sc1defJ} now becomes
\be \label{costD}
J(u)=\sum_{k=1}^N -x^k_2(t_1)\, x^k_3(t_1).
\ee

We define
\begin{align*}
f^0(x^k,u):=\begin{cases}
f^0_{\alpha}(x^k,u),\quad x^k\in X_{\alpha},\\
f^0_{\beta}(x^k,u),\quad x^k\in X_{\beta},
\end{cases}
\end{align*}
where
\begin{align*}
f^0_{\alpha}(x^k,u)&=-\bigl(a(x^k_2)+b(x^k_2)u+c_s x^k_2\bigr)x^k_3,\\
f^0_{\beta}(x^k,u)&=-\bigl(a(x^k_2)+b(x^k_2)u)x^k_3.
\end{align*}
Hence, to minimize \eqref{costD} is equivalent to minimize
\be\label{eqcost}
J(u)=\sum_{k=1}^N \int_{t_0}^{t_1}f^0(x^k,u)\, dt-\sum_{k=1}^N x^{k0}_2x^{k0}_3.
\ee

One of the goals of this section is to prove that there exists an optimal control for this optimal control
problem and that, if \eqref{t1>hatt0} and \eqref{Dsuse} hold,
every optimal control is bang-bang with only one switching time. More precisely, we prove the
following Theorem \ref{infimumJ} and Theorem \ref{application}.

Using \eqref{h(ys)>0}, we can easily check the continuity of the exit time with respect to 
the weak-$^*$ $L^{\infty}$ topology for the control. From 
the standard Arzel\`{a}-Ascoli 
theorem, we then get the following theorem (see also \cite[Theorem 1]{Piccoli1})
\begin{thm}\label{infimumJ}
The optimal control problem \eqref{sc1minJ} has a solution, i.e.,
there exists $u_*\in L^{\infty}((t_0,t_1);[w,1])$ such that
 \beq
 J(u_*)=\inf_{u\in L^{\infty}((t_0,t_1);[w,1])}J(u).
 \eeq
\end{thm}

\begin{thm}\label{application}
Let us assume that \eqref{t1>hatt0} and \eqref{Dsuse} hold. Then, for every optimal control
$u_*$ for the optimal control problem \eqref{sc1minJ},
there exists $t_*\in (t_0,t_1)$ such that
\begin{gather}
\label{bang-bang-effectif-2} u_*= w \text{ in } (t_0,t_*)\ \text {
and }\ u_*=1 \text{ in } (t_*,t_1).
\end{gather}
\end{thm}

This section is organized as follows. 
In subsection \ref{necessary} we prove a HMP (Theorem~\ref{PMP}) for our optimal control
problem. In subsection~\ref{bang}  we show how to deduce Theorem~\ref{application} from Theorem~\ref{PMP}.


\subsection{Hybrid Maximum Principle}\label{necessary}
Let us define the  Hamiltonian
\begin{equation*}
\begin{array}{ccc}
\mathcal{H}:  (\R^3)^N \times \R \times (\R^3)^N
&\rightarrow &\R
\\
(x,u,\psi)=((x^1,x^2,\ldots,x^N),u,(\psi^1,\psi^2,\ldots,\psi^N))& \mapsto & \mathcal{H}(x,u,\psi)
\end{array}
\end{equation*}
by
\begin{gather}
\label{ham0}
\mathcal{H}(x,u,\psi):=\sum_{k=1}^N \langle f(x^k,u),\psi^k \rangle
-\sum_{k=1}^N f^0(x^k,u).
\end{gather}
In \eqref{ham0} and in the following, $\langle a,b \rangle $ denotes the usual scalar product of $a\in \R^3$
and $b\in \R^3$. Let us also define the Hamilton-Pontryagin function $H:(\R^3)^N \times (\R^3)^N \rightarrow \R$ by
\begin{gather}
\label{ham02}
H(x,\psi):=\max_{u\in [w,1]} \mathcal{H}(x,u,\psi).
\end{gather}
It follows from \cite{Dmitruk, Piccoli2, TH, Piccoli1, Sussmann2, Sussmann1} that we have the following theorem:
\begin{thm}\label{PMP}
Let $u_*\in L^{\infty}((t_0,t_1);[w,1])$ be an optimal control for the optimal control problem
\eqref{sc1minJ}.
Let $x^k_*=(x^k_{*1},x^k_{*2},x^k_{*3})\tr$, $k=1,\cdots,N$, be the corresponding optimal trajectory, i.e.
$x^k_{*}\in\bigl( W^{1,\infty}(t_0,t_1)\bigr)^3$ are solutions to the following Cauchy problems
\begin{gather}
\label{eqdotx1}
{\dot x^k}_{*1}=1,\quad \quad x^k_{*1}(t_0)=x^{k0}_1,
\\
\label{eqdotx2}
{\dot x^k}_{*2}=a(x^k_{*2})+b(x^k_{*2})\, u_*,\quad x^k_{*2}(t_0)=x^{k0}_2,
\\
\label{eqdotx3}
{\dot x^k}_{*3}=c(x^k_{*2})\, x^k_{*3},\quad x^k_{*3}(t_0)=x^{k0}_3.
\end{gather}

If $y_s\in \left\{x^k_{*2}(t);\, t\in[t_0,t_1]\right\}$, let $\hat t_k\in[t_0,t_1]$
be the exit time for the $k$-th Dirac mass, i.e. the unique time $\hat t_k\in [t_0,t_1]$ such that
$x^k_{*2}(\hat t_k)=y_s$. 

If $y_s\notin\left\{x^k_{*2}(t);\, t\in (t_0,t_1]\right\}$, let $\hat t_k=t_1+1$.
Then, there exists $N$ vector functions $\psi^k=(\psi^k_1,\psi^k_2,\psi^k_3)\tr\in
\bigl(W^{1,\infty}(((t_0,\hat t_k)\cup(\hat t_k,t_1))\cap (t_0,t_1))\bigr)^3$,
such that
\begin{align}\label{N1}
& \dot\psi^k_1= 0,\\
\label{N2}
&\dot\psi^k_2=  -(a'(x^k_{*2})+b'(x^k_{*2}) u_*)\psi^k_2
-(a'(x^k_{*2})+b'(x^k_{*2}) u_*)x^k_{*3}\nonumber\\
&\qquad\qquad\qquad -c(x^k_{*2}) x^k_{*3}\quad \text{ in }\quad ((t_0,\hat t_k)\cup(\hat t_k,t_1))\cap (t_0,t_1),\\
\label{N3}
&\dot\psi^k_3= -c(x^k_{*2})\, \psi^k_3-(a(x^k_{*2})+b(x^k_{*2})u_*)- c(x^k_{*2})x^k_{*2},\\
\label{finalcondition-1-3}
&\psi^k_1(t_1)=\psi^k_3(t_1)=0,
\end{align}
and
\begin{align}
\label{psik1}
\psi^k_1(\hat t_k-0)&=\psi^k_1(\hat t_k+0),\\
\label{psik3}
\psi^k_3(\hat t_k-0)&=\psi^k_3(\hat t_k+0),
\end{align}
\begin{itemize}
\item if $\hat t_k<t_1$,
\begin{gather}
\label{use0N}
\psi^k_2(\hat t_k +0)-\psi^k_2(\hat t_k-0)\in  \Big[\f{c_s\, x^k_{*3}(\hat t_k)(y_s+\psi^k_3(\hat t_k))}{a(y_s)+b(y_s)},
\f{c_s\, x^k_{*3}(\hat t_k)(y_s+\psi^k_3(\hat t_k))}{a(y_s)+b(y_s)\, w}\Big],
\\
\label{finalcondition-2}
\psi^k_2(t_1)=0,
\end{gather}
\item if $\hat t_k=t_1$,
\begin{gather}
\label{finalcondition-2-caslimite}
-\psi^k_2(t_1)\in  \Big[0,
\f{c_s\, x^k_{*3}(t_1)y_s}{a(y_s)+b(y_s)\, w}\Big].
\end{gather}
\end{itemize}
Moreover, there exists a constant $h$ such that the following condition holds
\begin{gather}
\label{Hamiltonian}
\mathcal{H}(x^k_*(t),u_*(t),\psi^k(t))
=H(x^k_*(t), \psi^k(t))=h,\quad a.e.\quad t \in (t_0,t_1).
\end{gather}
\end{thm}

\textbf{Proof of Theorem \ref{PMP}.}
For sake of simplicity, 
we give the proof only for one Dirac mass $(N=1)$.
To simplify the notations we also
delete the $k = 1$ index.
For more than one Dirac mass, the proof is similar.

Applying the HMPs given in  \cite{Dmitruk, Piccoli2, TH, Piccoli1, 
Sussmann2, Sussmann1},
we get the existence of $\psi=(\psi_1,\psi_2,\psi_3)\tr\in
\bigl(W^{1,\infty}(((t_0,\hat t)\cup(\hat t,t_1))\cap 
(t_0,t_1))\bigr)^3$ such that  \eqref{N1} to \eqref{psik3} and, if $\hat 
t <t_1$, \eqref{finalcondition-2} hold, together with the existence of 
$h\in \R$ such that \eqref{Hamiltonian} is satisfied. Let us finally 
deal with \eqref{use0N} and \eqref{finalcondition-2-caslimite}. Let us 
treat only the case where $\hat t<t_1$ (the case $\hat t=t_1$ being 
similar). We follow \cite{TH}. From \eqref{Hamiltonian}, there exist 
$v_1\in [w,1]$ and $v_2\in [w,1]$ such that
\begin{align}
\label{H-}
H(\hat t-0)&=\max_{v\in[w,1]}\mathcal{H}(x(\hat t-0),v,\psi(\hat 
t-0)=H(x(\hat t-0),v_1,\psi(\hat t-0))\nonumber\\
&=<f(x(\hat t-0),v_1),\psi(\hat t-0)>-f^0(x(\hat t-0),v_1),\\
\label{H+}
H(\hat t+0)&=\max_{v\in[w,1]}\mathcal{H}(x(\hat t+0),v,\psi(\hat 
t+0))=H(x(\hat t+0),v_2,\psi(\hat t+0))\nonumber\\
&=<f(x(\hat t+0),v_2),\psi(\hat t+0)>-f^0(x(\hat t+0),v_2).
\end{align}
From \eqref{Hamiltonian}, \eqref{H-} and \eqref{H+}, we obtain
\be\label{eqq1}
\psi_2(\hat t+0)-\psi_2(\hat t-0)=\f{c_s x_3(\hat t\,) (y_s+\psi_3(\hat t\,))+b(y_s)(x_3(\hat t\,)+\psi_2(\hat t-0))(v_1-v_2)}{a(y_s)+b(y_s)v_2},
\ee
and
\be
\label{eqq2}
\psi_2(\hat t+0)-\psi_2(\hat t-0)=\f{c_s x_3(\hat t\,) (y_s+\psi_3(\hat t\,))+b(y_s)(x_3(\hat t\,)+\psi_2(\hat t+0))(v_1-v_2)}{a(y_s)+b(y_s)v_1}.
\ee
The Hamiltonian \eqref{ham0} becomes
\begin{align}
\label{calHN=1}
\mathcal{H}(x,u,\psi)=&(a(x_2)+c(x_2)x_2)\, x_3 +\psi_1 +a(x_2)\psi_2 +c(x_2)x_3\, \psi_3 \nonumber\\
&\qquad\qquad \qquad \qquad +b(x_2)(x_3+\psi_2)\, u,\quad t\in[t_0,t_1].
\end{align}
Let us denote 
\be\label{redefiPhi}
\Phi:=x_3+\psi_2.
\ee
When $t\neq\hat t$, from \eqref{scp}, \eqref{eqdotx3} and \eqref{N2},
we obtain
\be \label{tr3}
\f{d\Phi}{dt} =-(a'(x_2)+ b'(x_2)u)\, \Phi.
\ee
Noting that
\be\label{ruse}
\Phi(t_1)>0.
\ee
Combining \eqref{tr3} and \eqref{ruse}, we get
\be
\Phi(t)>0,\quad \forall t\in(\hat t,t_1].
\ee
Next, we analyze different cases:
\begin{enumerate}
\item
When $\Phi(\hat t-0)>0$ and $\Phi(\hat t+0)>0$, 
we have $v_1=v_2=1$. From \eqref{eqq1} or \eqref{eqq2}, we get
\be\label{b1}
\psi_2(\hat t+0)-\psi_2(\hat t-0)=\f{c_sx_3(\hat t)(y_s+\psi_3(\hat t))}{a(y_s)+b(y_s)}.
\ee
\item
When $\Phi(\hat t-0)<0$ and $\Phi(\hat t+0)>0$, 
we have $v_1=w$ and $v_2=1$. From \eqref{eqq1}, we get
\be\label{reo1}
\psi_2(\hat t+0)-\psi_2(\hat t-0)=\f{c_s x_3(\hat t)(y_s+\psi_3(\hat t))+b(y_s)\Phi(\hat t-0)(w-1)}{a(y_s)+b(y_s)}.
\ee
Since $\Phi(\hat t-0)<0$, from \eqref{reo1}, we obtain
\be\label{b3}
\psi_2(\hat t+0)-\psi_2(\hat t-0)>\f{c_s x_3(\hat t)(y_s+\psi_3(\hat t))}{a(y_s)+b(y_s)}.
\ee
From \eqref{eqq2}, we have
\be\label{reo2}
\psi_2(\hat t+0)-\psi_2(\hat t-0)=\f{c_s x_3(\hat t)(y_s+\psi_3(\hat t))+b(y_s)\Phi(\hat t+0)(w-1)}{a(y_s)+b(y_s)w}.
\ee
Since $\Phi(\hat t+0)>0$, from \eqref{reo2}, we obtain
\be\label{b4}
\psi_2(\hat t+0)-\psi_2(\hat t-0)<\f{c_s x_3(\hat t)(y_s+\psi_3(\hat t))}{a(y_s)+b(y_s)w}.
\ee
\item
When $\Phi(\hat t-0)=0$ and $\Phi(\hat t+0)>0$, from \eqref{eqq1}, we obtain
\be\label{b8}
\psi_2(\hat t+0)-\psi_2(\hat t-0)=\f{c_sx_3(\hat t)(y_s+\psi_3(\hat t))}{a(y_s)+b(y_s)}.
\ee
\end{enumerate}
In the three cases, we have proved that jump condition \eqref{use0N} holds.
This concludes the proof of Theorem \ref{PMP}.
\cqfd
\subsection{Proof of Theorem \ref{application}}\label{bang}
In this section, we use the necessary optimality conditions given in Theorem
\ref{PMP} to prove Theorem \ref{application}.
 From now on, we assume that the target time $t_1$ satisfies $t_1>\hat t_0$ so that all the cells will exit
 from Phase 1 into Phase 3
 before time $t_1$.
We give a proof of Theorem \ref{application} in the case where $N=1$ in section \ref{one}.
In section \ref{two}, we study the case where $N>1$; in this case
we need additionally to analyze the dynamics between different exit times $\hat t_k$, $k=1,2,\cdots,N$,
to obtain that there exists one and only one switching time and that the optimal switching direction is from $u=w$ to $u=1$.
In both cases $N=1$ or $N > 1$, we give
some numerical illustrations, respectively in section \ref{numericalpic1} and section \ref{numericalpic2}.
\subsubsection{Proof of Theorem \ref{application} in the case $N=1$}\label{one}
Let $u$ be an optimal control for the optimal control problem \eqref{sc1minJ} and let $x=(x_1,x_2,x_3)\tr$ be the corresponding trajectory.  Note that, by \eqref{ab}, $b(x_2)> 0$. Then, by \eqref{ham02}, \eqref{Hamiltonian}, \eqref{calHN=1} and \eqref{redefiPhi},
one has, for almost every $t\in (t_0,t_1)$,
\begin{gather}
\label{u=1} u(t)= 1\quad \text{ if }\quad \Phi(t)>0,
\\
\label{u=w} u(t)= w\quad \text{ if }\quad \Phi(t)<0.
\end{gather}
Let us recall that, under assumption \eqref{t1>hatt0} of Theorem~\ref{application}, there exists one and
 only one $\hat t \in [t_0,t_1)$ such that
\be
\label{x_2(hatt)=ys}
x_2(\hat t\, )=y_s.
\ee
Then
\be
\label{x2tt1}
x_2(t)>y_s, \quad  \forall t\in (\hat t, t_1].
\ee
We study the case where $\hat t >t_0$, the case $\hat t =t_0$ being obvious. Thanks to \eqref{use0N}, we get
\be \label{com1}
\Phi(\hat t+0)-\Phi(\hat t-0)\geqslant c_s\, x_3(\hat t\, )\, \f{y_s+\psi_3(\hat t\, )}{a(y_s)+b(y_s)}.
\ee
By \eqref{eqdotx2} and \eqref{N3}, we get
\be\label{rederive}
\f{d(x_2+\psi_3)}{dt}=-(x_2+\psi_3)\, c(x_2),
\ee
and then, using also \eqref{scp}, \eqref{finalcondition-1-3},  \eqref{x_2(hatt)=ys}, \eqref{x2tt1} and \eqref{rederive}, we obtain
\be \label{com2}
y_s+\psi_3(\hat t\, )=(x_2+\psi_3)(\hat t\, )=(x_2+ \psi_3)(t_1)\geqslant y_s.
\ee
Combining \eqref{com1} with \eqref{com2}, we get
\be\label{in1}
\Phi(\hat t-0)\leqslant \Phi(\hat t+0)-c_s\, x_3(\hat t\, )\f{y_s}{a(y_s)+b(y_s)}.
\ee
By \eqref{ab} and \eqref{tr3}, we obtain
\be \label{tr1}
\Phi(\hat t+0)=\Phi(t_1)\, e^{-\int_{\hat t}^{t_1}(2x_2(s)-c_1u(s))\, ds}.
\ee
Using the first inequality of \eqref{Dsuse}, \eqref{finalcondition-2}, \eqref{x2tt1} and \eqref{tr1}, we get
\be
\label{x3psi2}
\Phi(\hat t+0)\leqslant x_3(t_1).
\ee
Noticing  that $x_3(\hat t+0)=x_3(t_1)$ and using  \eqref{in1} and \eqref{x3psi2}, we get
\be
\label{tr2}
\Phi(\hat t-0)\leqslant x_3(t_1)(1-c_s\, \f{y_s}{a(y_s)+b(y_s)}).
\ee
 From the second inequality of \eqref{Dsuse} and
\eqref{tr2}, we get
\be
\Phi(\hat t-0)<0.
\ee
which, together with \eqref{tr3}, gives us
\be\label{tr4}
\Phi(t)<0, \quad t\in[t_0,\hat t\, ).
\ee
Moreover, by \eqref{finalcondition-2}, we have
\beq
\Phi(t)=(x_3+\psi_2)(t_1)=x_3(t_1)>0,
\eeq
which together with \eqref{tr3}, gives
\be\label{tr5}
\Phi(t)>0,\quad t\in(\hat t,t_1].
\ee
Taking $t_*=\hat t$ and combining \eqref{tr4} and \eqref{tr5}, with \eqref{u=1} and \eqref{u=w}, 
we conclude the proof of Theorem \ref{application} in the case where $N=1$.
\cqfd

\subsubsection{Numerical illustration in the case $N=1$}\label{numericalpic1}
\begin{figure}[htbp!]
\includegraphics[width=3in]{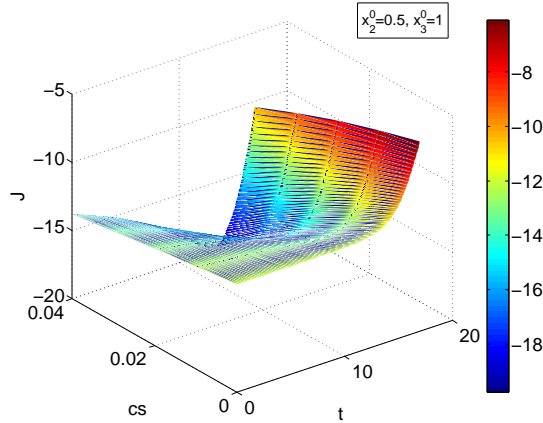}
\centering \vspace{-20mm}
\caption{Value of the cost function $J$ with respect to the switching time ($t$) and $c_s$ parameter in the case of one Dirac mass. When $c_s$ is ``small'', there is no switching time ($t=0$) and the optimal control is constant ($u=1$), while, when $c_s$ is ``large'', the optimal control strategy consists in switching from $u=w$ to $u=1$ at a time coinciding with the exit time. The initial values are specified in the insert.}
\label{continuous-fig}
\end{figure}

For one Dirac mass, the optimal switching time is unique.
Assumption \eqref{Dsuse} is not necessary to guarantee that the optimal control is a bang-bang control with only one
switching time.
It is just used to guarantee that the optimal switching time coincides with the exit time.
We give a numerical example to show that when $c_s$ is ``small", there is no switch at all and the optimal control is constant $(u=1)$,
while when $c_s$ is ``large", there is a switch occuring at the exit time (see figure~\ref{continuous-fig}).

The default parameter values are specified in Table \ref{tabparam} for the numerical studies. 
\setlength{\unitlength}{0.085in}

\begin{table}[ht]
\begin{center}
\begin{tabular}{| c | l | c |}
\hline
$t_0$&initial time &0.0\\
$t_1$&final time&17.0\\
$c_1$&slope in the $b(y)$ function&11.892\\
$c_2$&origin ordinate in the $b(y)$ function&2.288\\
$y_s$&threshold maturity&6.0\\
$w$&minimal bound of the control&0.5\\
\hline
\end{tabular}
\caption{Default parameter values}\label{tabparam}
\end{center}
\end{table}
\subsubsection{Proof of Theorem \ref{application} in the case $N>1$}\label{two}
Now, the Hamiltonian \eqref{ham0} becomes
\begin{align}
\label{calHN=N}
\mathcal{H}(x,u,\psi)&=\sum_{k=1}^N\bigl((a(x^k_{2})+c(x^k_2)x^k_{2})x^k_{3}+\psi^k_1
+a(x^k_2)\psi^k_2 +c(x^k_2)x^k_3\psi^k_3 \nonumber\\
&\qquad\qquad\qquad  \qquad+b(x^k_2)(x^k_3+\psi^k_2)u\bigr),\quad  t\in[t_0,t_1].
\end{align}
Reordering if necessary the $x^k$'s, we may assume, without loss of generality, that
\begin{gather}
\label{reordercroissant}
x^{10}_2<x^{20}_2<\ldots  < x^{k0}_2< \ldots <  x^{N0}_2.
\end{gather}
Let $u$ be an optimal control for the optimal control problem \eqref{sc1minJ} and let $x=(x^1, \ldots, x^k, \ldots x^N)$, with $x^k=(x^k_1,x^k_2,x^k_3)\tr$,  be the corresponding trajectory.
 From  \eqref{reordercroissant}, we have
\be
\label{hatdecroissant}
\hat t_N< \hat t_{N-1} < \ldots \hat t_k < \ldots < \hat t_1.
\ee
Let $\Phi_N :[t_0,t_1]\rightarrow \R$ be defined by
\begin{gather}
\label{Phi1} \Phi_N (t):=\sum_{k=1}^N b(x^k_2(t))(x^k_3(t)+\psi^k_2(t)).
\end{gather}
Noticing that $b(x^k_2)> 0$, by \eqref{ham02}, \eqref{Hamiltonian}, \eqref{calHN=N} and \eqref{Phi1}, one has, for almost every $t\in (t_0,t_1)$,
\begin{align}\label{optimalu}
u=w,\quad &\text{if}\quad \Phi_N (t)<0,\\
\label{optimalu2}
u=1,\quad &\text{if}\quad \Phi_N (t)>0.
\end{align}
We take the time-derivative of \eqref{Phi1} when $t\neq \hat t_k$, $k=1,\cdots,N$. From \eqref{ab}, we obtain
\begin{align}\label{dynamics1}
\dot \Phi_N  (t)=\sum_{k=1}^N (c_1(x^k_2)^2+2c_2x^k_2)(x^k_3+\psi^k_2).
\end{align}
Similarly to the above proof for one Dirac mass, we can prove that,
under assumption \eqref{Dsuse}, we have, for each $k=1,\cdots, N$, 
\begin{align}
\label{design1}
&(x^k_3+\psi^k_2)(t)<0,\quad \text{when}\ t\in(t_0,\hat t_k),\\
\label{design2}
&(x^k_3+\psi^k_2)(t)>0,\quad \text{when}\ t\in(\hat t_k,t_1).
\end{align}
By \eqref{hatdecroissant}, \eqref{Phi1}, \eqref{design1} and \eqref{design2}, and note that $b(x^k_2)>0$, we get
\begin{align}\label{DecideN1}
&\Phi_N (t)<0,\quad \text{when}\quad t\in(t_0,\hat t_N),\\
\label{DecideN2}
&\Phi_N (t)>0,\quad\text{when}\quad t\in(\hat t_1,t_1).
\end{align}
The key point now is to study the dynamics of $\Phi_N$ between different exit times $\hat t_k$.
Let $k\in \{1,\cdots,N-1\}$ and let us assume that
\be
\label{ua1}
\Phi_N (t)=0,\quad \text{for some}\quad t\in(\hat t_{k+1},\hat t_{k}).
\ee
 From \eqref{Phi1} and \eqref{ua1}, we get
\be
\label{x3+psi2}
x^k_3(t) +\psi^k_2(t)=-\sum_{i\neq k}\f{b(x^i_2(t))}{b(x^k_2(t))}(x^i_3(t) +\psi^i_2(t)).
\ee
 From \eqref{design1} and \eqref{design2}, for every $t\in(\hat t_{k+1},\hat t_k)$,
\begin{align}\label{decideN1}
&x^i_3(t)+\psi^i_2(t) <0, \quad \text{when}\quad i\leqslant k-1,\\
\label{decideN2}
&x^i_3(t)+\psi^i_2(t) >0, \quad \text{when}\quad i\geqslant k+1.
\end{align}
 From \eqref{ab}, \eqref{dynamics1} and \eqref{x3+psi2}, we get
\begin{align}\label{finalDN}
\dot\Phi_N (t)
=&\sum_{i\leqslant k-1}\f{x^i_3
+\psi^i_2}{b(x^k_2)}\bigl(c_1^2x^i_2x^k_2+2c_2^2+c_1c_2(x^i_2+ x^k_2)\bigr)(x^i_2-x^k_2)\nonumber\\
&+\sum_{i\geqslant k+1}\f{x^i_3
+\psi^i_2}{b(x^k_2)}\bigl(c_1^2x^i_2x^k_2+2c_2^2+c_1c_2(x^i_2+ x^k_2)\bigr)(x^i_2-x^k_2).
\end{align}
 From  \eqref{reordercroissant}, we get
\begin{align}\label{ua2}
&x^i_2(t)-x^k_2(t)<0,\quad \text{when}\quad i\leqslant k-1,\\
\label{ua3}
&x^i_2(t)-x^k_2(t)>0,\quad \text{when}\quad i\geqslant k+1.
\end{align}
Using  \eqref{decideN1} to \eqref{ua3}, we get
\be\label{finalN}
\dot\Phi_N (t)>0 \quad \text{whenever}\quad \Phi_N (t)=0, \quad \forall t\in(\hat t_{k+1},\hat t_k).
\ee
Combining \eqref{optimalu}, \eqref{optimalu2}, \eqref{DecideN1}, \eqref{DecideN2} and \eqref{finalN} together,  we get the existence of $t_*\in(t_0,t_1)$ such that
\beq
u_*= w \quad\text{in}\quad  (t_0,t_*) \quad \text{and}\quad  u_*= 1 \quad\text{in}\quad  (t_*,t_1).
\eeq
This concludes the proof of Theorem \ref{application}.
\cqfd

\subsubsection{Numerical illustration in the case $N>1$}\label{numericalpic2}

\begin{figure}[htbp!]
\setlength{\abovecaptionskip}{-1cm}
\begin{minipage}[t]{0.5\linewidth}
\centering
\includegraphics[width=2.2in]{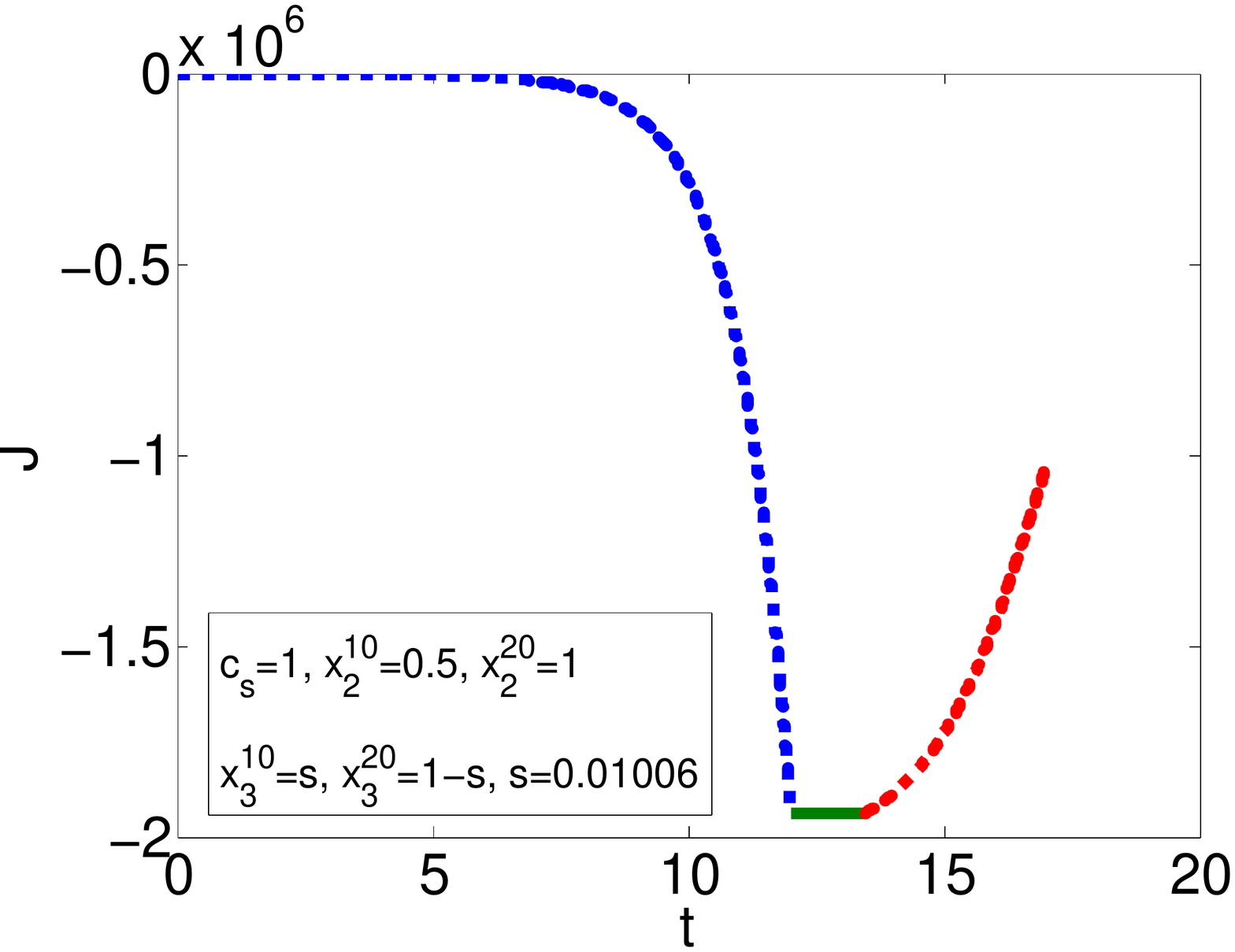}
\label{fig:side:a}
\end{minipage}%
\begin{minipage}[t]{0.5\linewidth}
\centering
\includegraphics[width=2.2in]{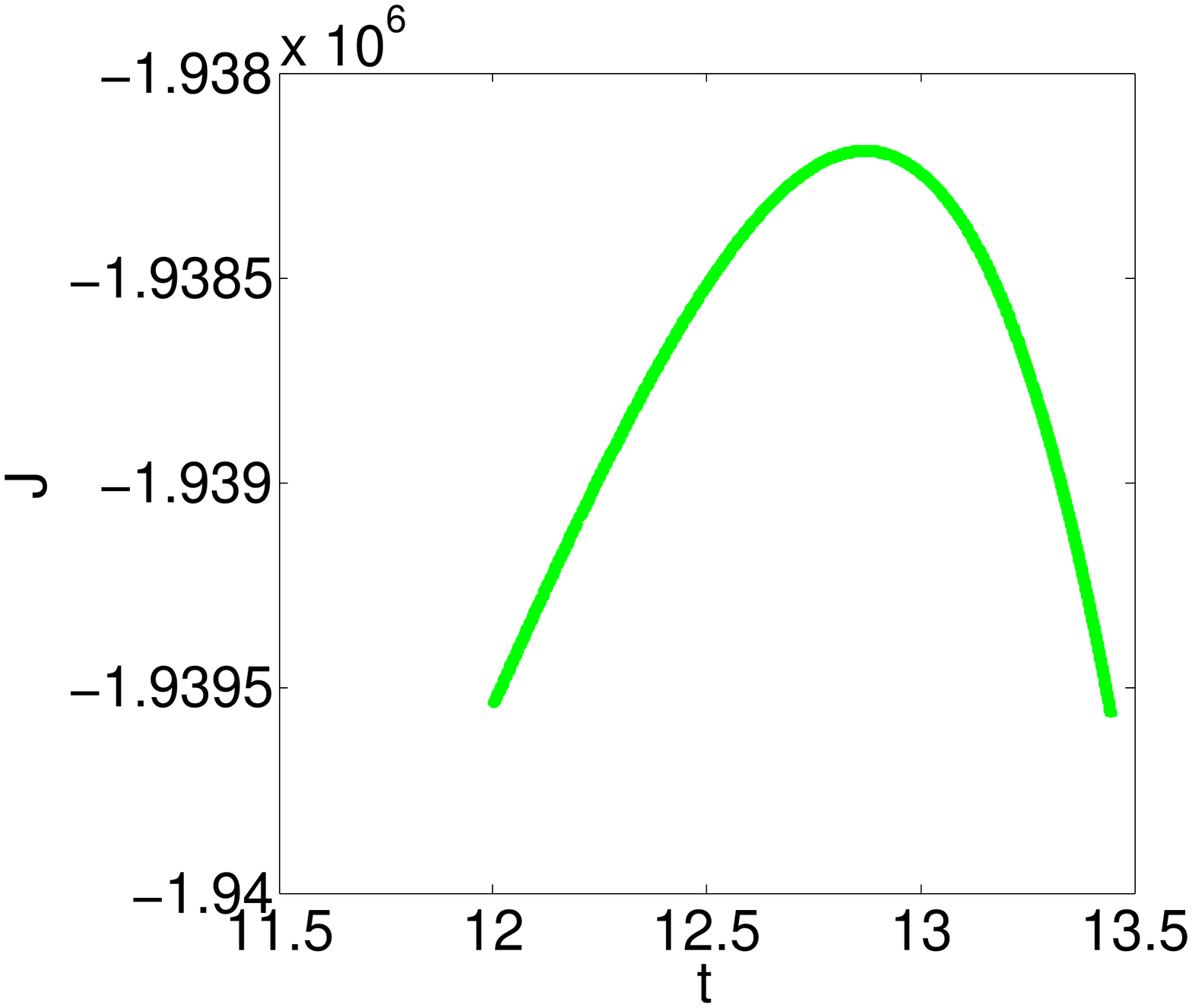}
\label{fig:side:b}
\end{minipage}
\caption{Value of the cost function $J$ with respect to the switching time ($t$) in the case of two Dirac masses and a ``large'' value of $c_s$ ($c_s=1.0$).
In the left panel,  the three-part curve represents the value of the cost function obtained after switching from $u=w$ to $u=1$ at time $t$. Blue dashed curve: switching time occurring before the first exit time;  green solid curve: switching time occurring in between the two exit times; red dashed curve: switching time occurring after the second exit time. The initial values are specified in the insert.
The right panel is a zoom on the green solid curve displayed on the left panel. There are two optimal switching times which coincide with the two exit times.}
\label{otw}
\end{figure}

\begin{figure}[htbp!]
\setlength{\abovecaptionskip}{-1cm}
\begin{minipage}[t]{0.5\linewidth}
\centering
\includegraphics[width=2.2in]{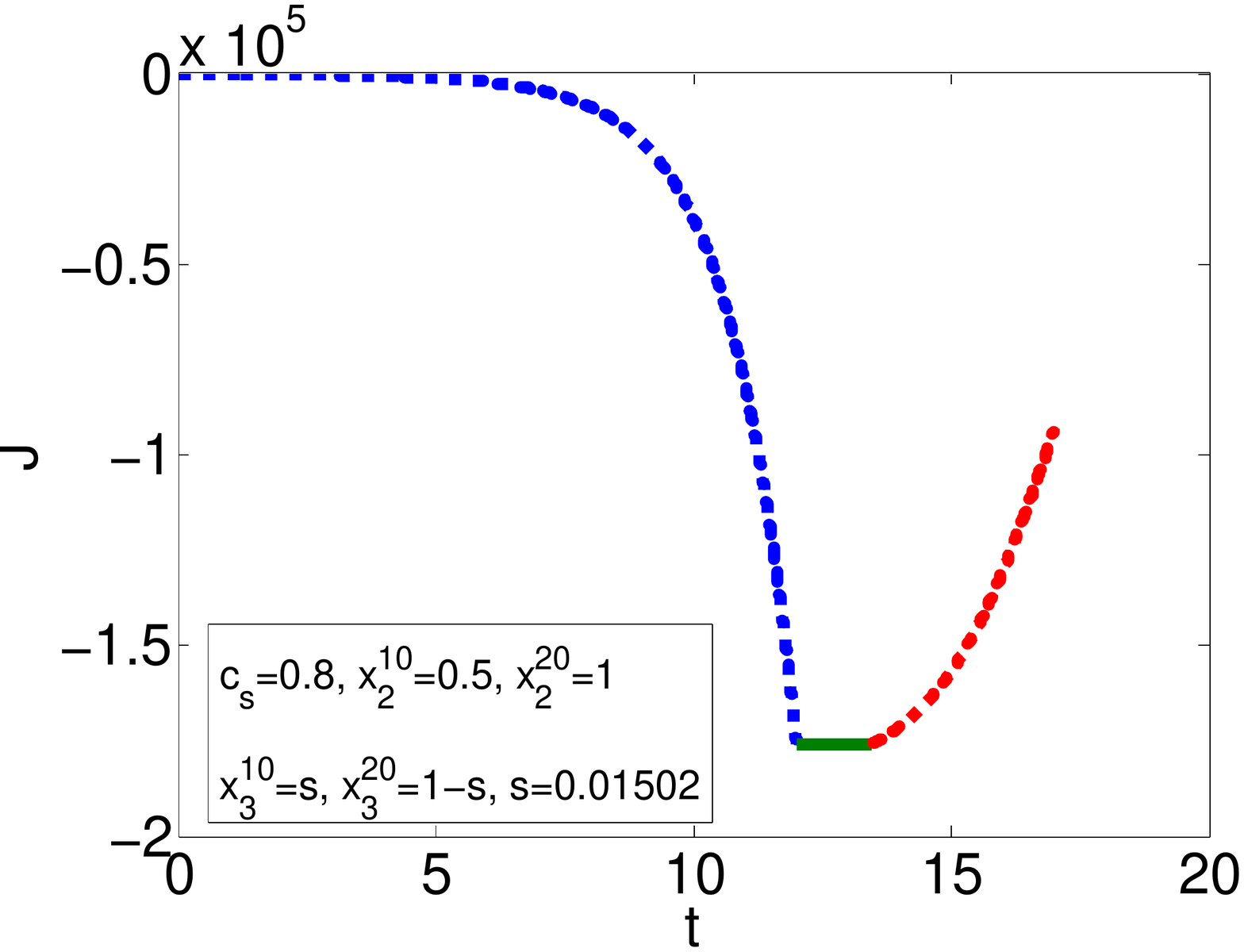}
\label{fig:side:a-cs-large}
\end{minipage}%
\begin{minipage}[t]{0.5\linewidth}
\centering
\includegraphics[width=2.2in]{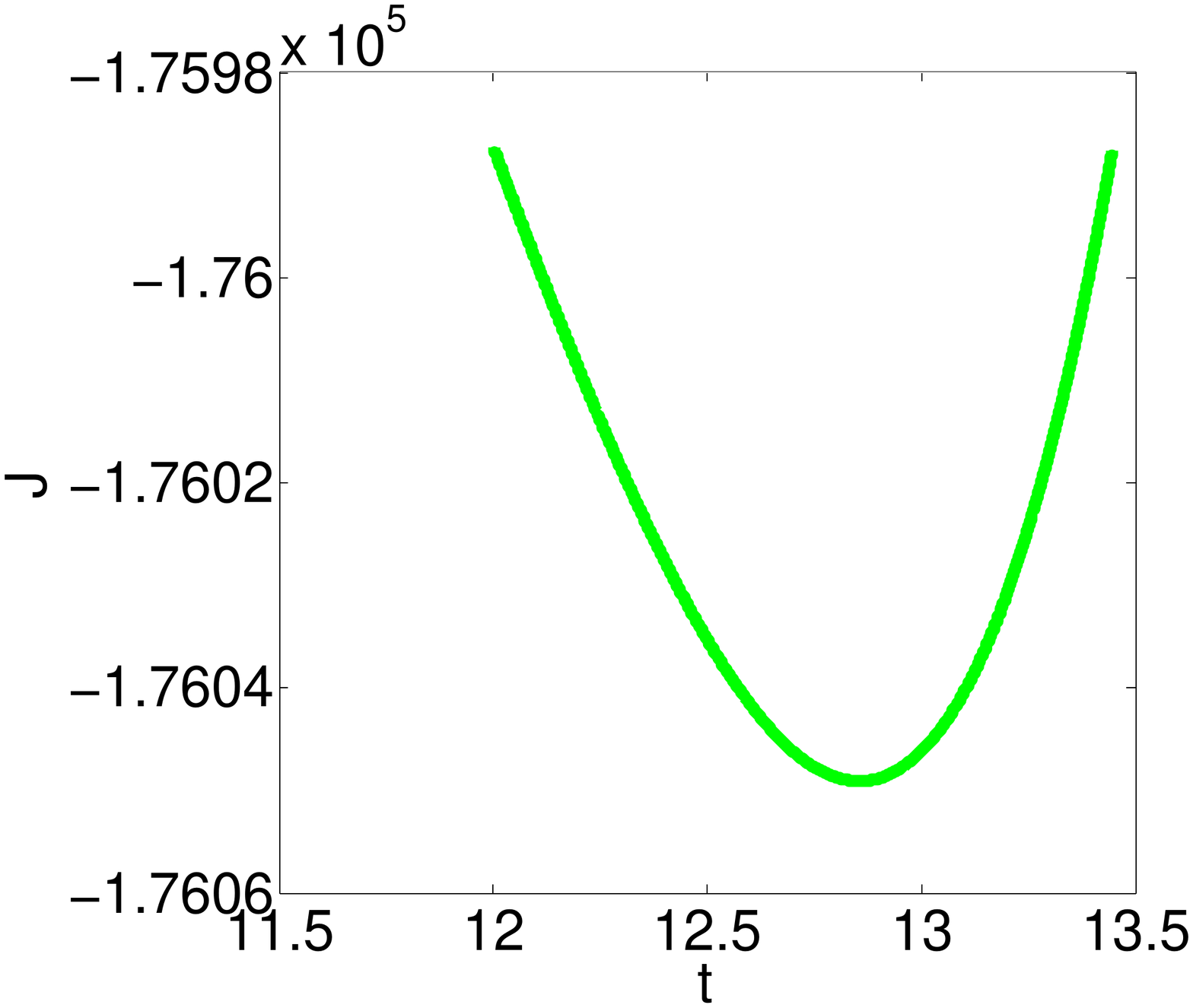}
\label{fig:side:b-cs-large}
\end{minipage}
\caption{Value of the cost function $J$ with respect to the switching time ($t$) in the case of two Dirac masses and a ``large'' value of $c_s$ ($c_s =0.8$).
In the left panel,  the three-part curve represents the value of the cost function obtained after switching from $u=w$ to $u=1$ at time $t$. Blue dashed curve: switching time occurring before the first exit time; green solid curve: switching time occurring in between the two exit times;  red dashed curve: switching time occurring after the second exit time. The initial values are specified in the insert.
The right panel is a zoom on  the green solid curve displayed on the left panel. There is one single optimal switching time, which occurs in between the two exit times.}
\label{btw}
\end{figure}

The optimal control can in some cases be not unique for more than one Dirac mass.
Let us consider the case of two Dirac masses as an example.
The optimal switching time may happen either at the first exit time or at the second exit time (see figure~\ref{otw}),
or between the two exit times (see figure~\ref{btw}).

\begin{rem}
The values of the cost function $J$ in Fig \ref{continuous-fig}
and Fig \ref{otw}, Fig \ref{btw} have different orders.
This is due to differences in the value of the proliferation rate $c_s$.
There is a great contribution of the cell mass to the criterion when $c_s$ is ``large" in Fig \ref{otw} and Fig \ref{btw}.
\end{rem}
\section{Optimal control in the PDE case}\label{Pde}
In this section, we study the optimal control in the PDE case.
We give the proof of Theorem~\ref{Theorempde}.
We first give an explicit expression for the cost function $J$ defined in \eqref{sc1defJ}.

Let us  define a new map
\begin{equation*}
\begin{array}{ccc}
e: [0,y_s]\times L^\infty ((t_0,t_1);[w,1])
&
\rightarrow
&
[t_0,t_1]
\\
(y_0,u)
&
\mapsto
&
e (y_0,u)
\end{array}
\end{equation*}
by requiring $\Psi(e(y_0,u),y_0,u)=y_s$,
where $\Psi$ is defined by \eqref{defPhi}.
Note that, under assumption \eqref{t1>hatt0}, one has, for every $y_0\in [0,y_s]$, the existence of $t\in [t_0,t_1]$
such that
\begin{gather}
\Psi(t,y_0,u)=y_s.
\label{defg-existence}
\end{gather}
Again, \eqref{h(ys)>0} implies that there exists at most one $t\in[t_0,t_1]$ such that
\eqref{defg-existence} holds. This shows that $e$ is well defined. Moreover, we have the following lemma
\begin{lem}
\label{lemmacontinuityg}
Let $(y^n_0)_{n\in \mathbb{N}}$ be a sequence of elements in $[0,y_s]$ and $(u^n)_{n\in \mathbb{N}}$ be a sequence of elements in $L^\infty ((t_0,t_1);[w,1])$. Let us assume that, for some $y_0\in [0,y_s]$
and for some $u\in L^\infty ((t_0,t_1);[w,1])$,
\begin{gather*}
y^n_0\rightarrow y_0 \text{ as } n\rightarrow +\infty,
\\
u^{n}\overset{*}{\rightharpoonup} u \text{ in } L^{\infty}(t_0,t_1) \text{ as } n\rightarrow +\infty.
\end{gather*}
Then
\beq
e(y^n_0,u^{n})\rightarrow e(y_0,u) \text{ as } n\rightarrow +\infty.
\eeq
\end{lem}

Let now $\rho_0$ be a Borel measure on $\R\times \R $ such that \eqref{rho0>0} and \eqref{supportrho0} hold.
Using \eqref{solution},  \eqref{sc1defJ} becomes
\be
\label{sc2}
J(u)=-\iint_{[0,1]\times[0,y_s]} \Psi(t_1,y_0,u)\,  e^{c_s e(y_0,u)}\,
d\rho_0(x_0,y_0).
\ee
In order to emphasize the dependence of $J$ on the initial data $\rho_0$, from now on we write
$J(\rho_0,u)$ for $J(u)$.

It is well known that there exists a sequence
$((x^{i,n}_0,y^{i,n}_0,\lambda^{i,n}_0))_{1\leqslant i\leqslant n,\, n\in \mathbb{N}}$ of elements
in $[0,1]\times [0,y_s]\times (0,+\infty)$ such that, if
\begin{gather}
\rho^n_0 :=\sum_{i=1}^n \lambda^{i,n}_0 \delta_{x^{i,n}_0,y^{i,n}_0},
\end{gather}
then
\begin{multline}
\label{limDirac}
\lim_{n\rightarrow +\infty}\iint_{[0,1]\times[0,y_s]}\varphi(x_0,y_0) \; d\rho^n_0(x_0,y_0) =
\\
\iint_{[0,1]\times[0,y_s]}\varphi(x_0,y_0)  d\rho_0(x_0,y_0),
\, \forall \varphi\in C^0([0,1]\times [0,y_s]).
\end{multline}
 From Theorem~\ref{infimumJ} and Theorem~\ref{application}, there exists $t^n_*\in [t_0,t_1]$ such that, if $u^n_*:[t_0,t_1]\rightarrow  [w,1]$ is
defined by
\begin{gather}
\label{defun*=}
u^n_*=w \text{ in } [t_0,t^n_*) \text{ and } u^n_*=1 \text{ in } (t^n_*,t_1],
\end{gather}
then
\begin{gather}
\label{Jun*}
J(\rho_0^n,u^n_*)\leqslant J(\rho_0^n,u),\quad \forall u\in  L^\infty ((t_0,t_1);[w,1]).
\end{gather}
Extracting a subsequence if necessary,
we may assume without loss of generality the existence of $t_*\in [t_0,t_1]$ such that
\begin{gather}
\label{cvtn*tot*}
\lim_{n\rightarrow +\infty}t^n_*=t_*.
\end{gather}
Let us define  $u_*:[t_0,t_1]\rightarrow  [w,1]$  by
\begin{gather}
\label{defu*=}
u_*=w \text{ in } [t_0,t_*) \text{ and } u_*=1 \text{ in } (t_*,t_1].
\end{gather}
Then, using  \eqref{defun*=}, \eqref{cvtn*tot*} and \eqref{defu*=}, one gets
\begin{gather}
\label{cvuniformPhi}
\Psi(t_1,\cdot,u^n_*)\rightarrow \Psi(t_1,\cdot,u_*) \text{ in }
C^0([0,y_s]) \text{ as } n\rightarrow +\infty.
\end{gather}
Moreover, from \eqref{defun*=}, \eqref{cvtn*tot*} and \eqref{defu*=}, one has
\begin{gather}
\label{cvweakun}
u^{n}_*\overset{*}{\rightharpoonup} u_* \text{ in } L^{\infty}(t_0,t_1) \text{ as } n\rightarrow +\infty.
\end{gather}
 From Lemma~\ref{lemmacontinuityg} and \eqref{cvweakun}, one gets
\begin{gather}
\label{cvgn}
e(\cdot,u^n_*)\rightarrow e(\cdot,u_*) \text{ in } C^0([0,y_s]) \text{ as } n\rightarrow +\infty.
\end{gather}
 From \eqref{sc2}, \eqref{limDirac}, \eqref{cvuniformPhi} and \eqref{cvgn} and a classical theorem on the weak topology (see, e.g.,
\cite[(iv) of Proposition 3.13, p. 63]{2011-Brezis-book}), one has
\begin{gather}
\label{cvJun}
J(\rho^n_0,u^n_*)\rightarrow J(\rho_0,u_*)  \text{ as } n\rightarrow +\infty.
\end{gather}
Let now $u\in  L^\infty ((t_0,t_1);[w,1])$.
 From Lemma~\ref{lemmacontinuityg}, \eqref{sc2} and \eqref{limDirac}, one gets
\begin{gather}
\label{limrhonu}
J(\rho^n_0,u)\rightarrow J(\rho_0,u)  \text{ as } n\rightarrow +\infty.
\end{gather}
Finally, letting $n\rightarrow +\infty$ in  \eqref{Jun*} and using \eqref{cvJun}
together with \eqref{limrhonu}, one has
$$
J(\rho_0,u_*)\leqslant J(\rho_0,u),
$$
which concludes the proof of Theorem \ref{Theorempde}.
\cqfd

\section*{Appendix}
\textbf{Sketch of another proof of Theorem \ref{PMP}}

In this section, we sketch another proof of Theorem \ref{PMP}, using approximation arguments inspired from \cite{S}.
The interest of this approach is that it might be more suitable to prove a maximal principle also in the PDE case.
For sake of simplicity, 
we show the proof only for one Dirac mass.
The idea is first to construct a smooth optimal control problem.
For the smooth optimal control problem, we can apply PMP.
By passing to the limit, we
then derive necessary optimality conditions for our discontinuous problem.

\textbf{Step 1.}
Let us denote by $\chi: \R \rightarrow \R$ the characteristic function of $(-\infty,y_s)$, i.e.
\be
\label{defchi}
\chi(x_2)=\begin{cases}
1,\quad \forall x_2\in (-\infty,y_s),\\
0,\quad \forall x_2\in [y_s,+\infty).
\end{cases}
\ee
Let $(w_i)_{i\in \mathbb{N}^*}$ be a sequence  of elements
in  $C^{\infty}(\R)$ such that
\be
\label{barwi-1}
0\leqslant  w_i,\quad \int_{\R} w_i(x)\, dx=1,\quad \text{support}\, w_i\subset [-1/i,0],\quad \forall i \in \mathbb{N}^*,
\ee
and, for some $C>0$,
\be
\label{barwi-2}
|w'_i(x)|\leqslant Ci^2, \quad \forall x \in \R, \quad \forall i\in \mathbb{N}^*,
\ee
(clearly such a sequence does exist).
Then, we define a sequence of functions $(\chi_i)_{i\in \mathbb{N}^*}$ from $\R$ into $\R$ as follows:
\be \label{barchii}
\chi_i(x):=\int_{\R}\chi(y)\, w_i(x-y)\, dy= \int_{-\infty}^{y_s}w_i(x-y)\, dy=
\int_{x-y_s}^{+\infty}w_i(z)\, dz, \, \forall i\in \mathbb{N}^*, \, \forall x\in \R.
\ee
Let $f_i:\R^3\times \R\rightarrow \R^3$ be defined by
\be
\label{deffi}
f_i(x,u):=\left(
\begin{array}{c}
 1  \\
 a(x_2)+b(x_2)\, u \\
 c_s\chi_i(x_2)\, x_3
\end{array}
\right),
\quad \forall x=(x_1,x_2,x_3)\tr \in \R^3,\quad \forall u\in \R.
\ee

Let $u_*$ be an optimal control for the optimal control
problem \eqref{sc1minJ} and let $x_*$ be the associated trajectory.
Let $(z_i)_{i\in \mathbb{N}^*}$ be  a sequence of uniformly bounded elements of $C^1([t_0,t_1])$ such that
\begin{align}
z_i\rightarrow u_*\ \text{in}\ L^2(t_0,t_1)\ \text{as}\ i\rightarrow + \infty.
\end{align}
Let us then define $J_i:L^{\infty}((t_0,t_1);[w,1])\rightarrow \R$ by
\beq
J_i(u):=\displaystyle-\int_{t_0}^{t_1}\bigl(a(x_2)+b(x_2)u+c_s\chi_i(x_2)x_2\bigr)x_3\, dt+
\f{1}{\sqrt{i}}\int_{t_0}^{t_1}|u(t)-z_i(t)|^2\, dt -x^0_2x^0_3,
\eeq
where $x:[t_0,t_1]\rightarrow \R^3$ is the solution to the Cauchy problem
\begin{gather}\label{dyxi}
\dot x=f_i(x,u),\,x(t_0)=x^0.
\end{gather}
We consider the following optimal control problem
\begin{equation*}
\tag{$\mathcal{P}_i$}
\text{ minimize } J_i (u) \text{ for } u \in L^{\infty}((t_0,t_1);[w,1]).
\end{equation*}
 For any $i=1,2,\cdots$, problem $(\mathcal{P}_i)$ is a ``smooth'' optimal control problem.
By a classical result in optimal control theory
(see, e.g., \cite[Corollary 2, p. 262]{LM}), there exists an optimal control $u_i$ for problem $(\mathcal{P}_i)$.
Let $x_i$ be the optimal trajectory corresponding to the control $u_i$ for dynamics \eqref{dyxi}. We have the
following lemma (compare to \cite[Lemma 4]{S}):
\begin{lem}\label{lem9}
The following holds as $i\rightarrow + \infty$
\begin{gather}
\label{condition1}
u_i\rightarrow u_* \text{ in } L^2(t_0,t_1),
\\
\label{condition2}
x_i\rightarrow x_* \text{ in } C^0([t_0,t_1];\R^3),
\\
\label{chixi2}
\chi_i(x_{i2})\rightarrow \chi(x_{*2})\text{ in } L^1(t_0,t_1).
\end{gather}
\end{lem}
\textbf{Step 2.}  We now deduce necessary optimality conditions for the optimal
control problem \eqref{sc1minJ} in the form of PMP.
The Hamiltonian and the Hamilton-Pontryagin function  for problem $(\mathcal{P}_i)$ are respectively
\begin{gather}
\label{defcalHi}
\mathcal{H}_i(x,u,\psi)=\langle f_i(x,u),\psi \rangle +(a(x_2)+b(x_2)u+c_s\chi_i(x_2)x_2)x_3
-\f{1}{\sqrt{i}}|u-z_i(t)|^2,
\\
\label{defHi}
H_i(x,\psi)=\max_{u\in [w,1]}\mathcal{H}_i(x,u,\psi).
\end{gather}
By the PMP -see, e.g., \cite[Theorem 2, p. 319]{LM} or \cite[Section 6.5]{BP}-, there exists an
absolutely continuous function $\psi_i:[t_0,t_1]\rightarrow \R^3$ such that
\begin{align}
\dot\psi_i\stackrel{a.e.}{=}&-\Big[\f{\partial f_i}{\partial x}(x_i,u_i)\Big]\tr \psi_i
-\f{\partial }{\partial x}\bigl((a(x_{i2})+b(x_{i2})u_i+c_s\chi_i(x_{i2})x_{i2})x_{i3}\bigr),
\label{necessary1}\\
\psi_i(t_1)= & \ 0,\label{necessary2}
\end{align}
and there exist constants $h_i$ such that
\be\label{maximum}
\mathcal{H}_i(x_i(t),u_i(t),\psi_i(t)) = H_i(x_i(t),\psi_i(t))=h_i,\quad a.e.\quad t \in (t_0,t_1).
\ee
Let us denote $\psi_i=(\psi_{i1},\psi_{i2},\psi_{i3})\tr$.
 From \eqref{deffi}, \eqref{necessary1} and \eqref{necessary2}, we have
\begin{align}\label{ew1}
&\dot\psi_{i1}=0,\\
\label{ew2}
&\dot\psi_{i2}=-(a'(x_{i2})+b'(x_{i2})\, u_i)\, \psi_{i2}-c_s\, \chi_i'(x_{i2})\, x_{i3}\,\psi_{i3}
-(a'(x_{i2})+b'(x_{i2})\, u_i)\, x_{i3}\nonumber\\
&\qquad\qquad\qquad\qquad\qquad \qquad -c_s\, \chi_i(x_{i2})\, x_{i3}-c_s\, x_{i2}\chi_i'(x_{i2})\, x_{i3},\\
\label{ew3}
&\dot\psi_{i3}=-c_s\, \chi_i(x_{i2})\,\psi_{i3}-(a(x_{i2})+b(x_{i2})u_i)- c_s\, x_{i2}\, \chi_i(x_{i2}),\\
\label{conditioniniti}
&\psi_{i1}(t_1)=\psi_{i2}(t_1)=\psi_{i3}(t_1)=0.
\end{align}
We can prove that
\be
\label{psii1=psi1}
\psi_{i1}(t)=\psi_1(t)=0,\quad \forall t\in[t_0,t_1],
\ee
and
\be\label{cvpsi3}
\psi_{i3} \rightarrow \psi_3 \text{ in } C^0([t_0,t_1]) \text{ as } i\rightarrow +\infty.
\ee
As far as $\psi_{i2}$ is concerned, 
Theorem~\ref{PMP}  in
the case where $x_{*2}(t_0)=x^0_2=y_s$ or $x_{*2}(t_1)< y_s$ follows directly from the standard PMP. Hence, we may assume that
\begin{gather}
\label{hypexit-new-large}
 x_{*2}(t_0)<y_s \leqslant x_{*2}(t_1).
\end{gather}
Let us treat the case where
\begin{gather}
\label{hypexit-new}
 x_{*2}(t_0)<y_s<x_{*2}(t_1),
\end{gather}
(the case
$x_{*2}(t_1)=y_s$
being similar).
By \eqref{h(ys)>0}, 
there exists one and only one  $\hat t\in (t_0,t_1)$ such that
\be
\label{cm0-1}
x_{*2}(\hat t\, )=y_s.
\ee
Using \eqref{condition2} and \eqref{cm0-1}, one also gets that, at least if $i$ is large enough, which, from now on, will always  be  assumed,
there exists one and only one  $\hat t_i\in (t_0,t_1)$ and one and only one $\bar t_i\in (t_0,t_1)$ such that
\be
\label{cm0-2}
x_{i2}(\hat t_i)=y_s,\quad
x_{i2}(\bar t_i)=y_s-(1/i).
\ee
Using \eqref{condition1} and \eqref{condition2}, we can prove
\begin{gather}
\label{hattibarti}
\hat t_i \rightarrow \hat t \text{ and } \bar t_i \rightarrow \hat t \text{ as } i\rightarrow +\infty.
\end{gather}
It is easy to check that
\be
\label{cm9}
\psi_{i2}\rightarrow \psi_2\quad \text{in}\quad C^0([t_0,\hat t-\varepsilon]\cup [\hat t+\varepsilon,t_1]),
\quad \forall \varepsilon>0.
\ee
We now prove jump condition \eqref{use0N} when $\hat t<t_1$, the 
proof of \eqref{finalcondition-2-caslimite} when $\hat t=t_1$ being similar.
Let us integrate \eqref{ew2} from $\bar t_i$ to $\hat t_i$,
 we get
\begin{gather}
\label{eqpsi-saut-1}
\psi_{i2}(\hat t_i)-\psi_{i2}(\bar t_i)=A(i) + B(i),
\end{gather}
with
\begin{gather}
\label{partiesanssaut}
A(i):=-\int^{\hat t_i}_{\bar t_i}
((a'(x_{i2})+b'(x_{i2})\, u_i)\, (\psi_{i2}+ x_{i3})+c_s\, \chi_i(x_{i2})\, x_{i3})\, dt,
\\
\label{partiessaut}
B(i):=-\,\int^{\hat t_i}_{\bar t_i}
c_s \, x_{i3}\,(x_{i2}+\psi_{i3})\,\chi_i'(x_{i2})\, dt.
\end{gather}
It is easy to obtain that
\begin{gather}
\label{estA}
A(i)\rightarrow 0 \text{ as } i\rightarrow +\infty.
\end{gather}
 For $B(i)$, we perform the change of variable $\tau=x_{i2}(t)$. By \eqref{cm0-2} and \eqref{partiessaut}, we get
\begin{gather}
\label{newexpBi}
B(i)=-\,\int^{y_s}_{y_s-(1/i)}
\frac{c_s \, x_{i3}(x_{i2}^{-1}(\tau))\, (\tau +
\psi_{i3}(x_{i2}^{-1}(\tau)))}{a(\tau)+b(\tau)u(x_{i2}^{-1}(\tau))}\,\chi_i'(\tau)\, d\tau.
\end{gather}
Let us point out that, from \eqref{barwi-1} and \eqref{barchii}, one has
\begin{gather}
\label{intchiprime}
\int^{y_s}_{y_s-(1/i)} \chi_i'(\tau)\, d\tau =-1, \quad  \chi_i'\leqslant 0.
\end{gather}
 From \eqref{h(ys)>0}, \eqref{condition2}, \eqref{cvpsi3},
  \eqref{newexpBi}, \eqref{intchiprime}, one gets that
\begin{gather*}
\frac{c_s\, x_{*3}(\hat t\, )(y_s+\psi_3(\hat t\, ))}{a(y_s)+b(y_s)} \leqslant \liminf_{i\rightarrow +\infty} B(i)\leqslant \limsup_{i\rightarrow +\infty} B(i)\leqslant \frac{c_s\, x_{*3}(\hat t\, )(y_s+\psi_3(\hat t\, ))}{a(y_s)+b(y_s)w},
\end{gather*}
which, together with \eqref{hattibarti}, \eqref{cm9}, \eqref{eqpsi-saut-1} and \eqref{estA}, gives \eqref{use0N}.

Letting $i\rightarrow +\infty$ in \eqref{maximum}, we get the 
existence of $h$ such that \eqref{Hamiltonian} holds. This concludes the proof of Theorem \ref{PMP}.
\cqfd

\section*{Acknowledgements}
We thank Emmanuel Tr\'{e}lat for useful discussions on the Hybrid Maximum Principle.

\end{document}